\documentclass{amsart}
\usepackage{amssymb}
\usepackage{amsmath}
\usepackage{amsfonts}
\usepackage{amsthm}

\newtheorem{theorem}{Theorem}[section]

\newtheorem{proposition}[theorem]{Proposition}
\newtheorem{corollary}[theorem]{Corollary}
\newtheorem{definition}[theorem]{Definition}
\newtheorem{example}[theorem]{Example}

\def\rmark{\mbox{$\rm\bf\rule{0.06em}{1.45ex}\kern-0.05em R$}}

\def\b1{b^{-1}}
\def\a1{a^{-1}}

\begin{document}

\begin{center}
{\textbf{\Large An Introduction to Bipolar Fuzzy Soft Hypervector Spaces}}

\vspace{5mm}

{\textbf{O. R. Dehghan}}\\[0pt]
\vspace{0.3cm}

\begin{minipage}{12.5cm}
\vspace{0.5cm} \textbf{ Abstract.} {\small The aim of this paper is to introduce the notion of bipolar fuzzy soft hypervector spaces and study their basic properties. In this regard, at first some new operation and external hyperoperation are defined on bipolar fuzzy soft sets over hypervector space $V$, related to the operation and external hyperoperation of $V$. Then the notion of bipolar fuzzy soft hypervector space is defined, supported by non-trivial examples, and it is investigated that the new bipolar fuzzy soft sets, constructed by the mentioned operation and hyperoperation, are bipolar fuzzy soft hypervector spaces. Finally, the behavior of bipolar fuzzy soft hypervector spaces under linear transformations is studied.}
\end{minipage}

\vspace{0.5cm}

\end{center}

\section{\protect\Large Introduction}
Fuzzy set, was introduced by Zadeh \cite{Zadeh} in 1965, is a mathematical tool for representing objects whose boundary is vague. In fuzzy sets, membership degrees indicate the degree of belonging of elements to the set or the degree of satisfaction of elements to the property related to the set.

There are several extensions for fuzzy sets. Fuzzy sets of type 2 represent membership degrees with fuzzy sets. $L$-fuzzy sets enlarge the range of membership degree $[0,1]$ into a lattice $L$. Interval-valued fuzzy sets represent the membership degree with interval values to reflect the uncertainty. Intuitionistic fuzzy sets, was first introduced by Atanassov \cite{Atanassov}, have a pair of membership degree and nonmembership degree. Each of these concepts, while having similarities and differences with the others, has its own uses.

The concept of bipolar fuzzy sets was introduced by Zhang \cite{Zhang} in 1994 for cognitive modeling and multi-agent decision analysis, as an extension of fuzzy sets, where the image of membership function is increased from $[0,1]$ to $[-1,1]$. Bipolar fuzzy sets have a pair of membership degrees that represent the degree of satisfaction to the property related to a fuzzy set and its counter-property. This extension has similarities and differences with the previous ones in semantics and representations. In this type of fuzzy sets, the membership degree $0$ means that elements are irrelevant to the corresponding property, the membership degrees on $(0,1]$ indicate that elements somewhat satisfy the property, and the membership degrees on $[-1,0)$ indicate that elements somewhat satisfy the implicit counter-property.

Also, the theory of soft sets, as an another mathematical tool for modeling uncertainty, was introduced by Molodtsov \cite{Molodtsov} in 1999. This theory was studied in various branches, particularly in algebraic structures; for examples, its application in a decision making problem by Maji \cite{Maji,Maji Roy}, soft groups by Aktas \cite{Aktas}, soft rings by Acar \cite{Acar} and soft vector spaces by sezgin \cite{Sezgin}.

Moreover, Cogman \cite{Cogman} in 2011 defined a fuzzy soft set, which is a more accurate tool for modeling, and studied its related properties. Then the theory of bipolar fuzzy soft sets was introduced and studied as a new applied generalization of previous theories. In fact, Abdollah \cite{Abdullah} in 2014 combined the concepts of bipolar fuzzy set and soft set and gave an general algorithm to solve decision making problems by using bipolar fuzzy soft sets. This idea has attracted the attention of some researchers and so many papers have been published based on it. For example, Akram \cite{Akram,Akram 2} studied applications of bipolar fuzzy soft sets in K-algebras and presented two techniques for the diagnosis of a disease in which the pairwise comparison of diseases and symptoms is considered in bipolar behavior, Ali \cite{Ali} used the technique of parameter reductions of bipolar fuzzy soft sets to solve decision-making problems, Abughazalah \cite{Abughazalah} applied the bipolar fuzzy sets in BCI-algebras, Riaz \cite{Riaz} discussed bipolar fuzzy soft topology with decision-making, Mahmood \cite{Mahmood} introduced the notion of bipolar complex fuzzy soft set as a generalization of bipolar complex fuzzy set and soft set and recently, Khan \cite{Khan} has applied the bipolar fuzzy soft matrices for solving the decision making problems.

On the other hand, algebraic hyperstructures was born in 1934, when Marty \cite{Marty} generalized the notion of operation into hyperoperation. An operation assigns to any two elements of the context set a unique element of that set, while a hyperoperation assigns to any two elements of the set a unique subset of that set. The theory of algebraic hyperstructures has been studied by many researchers in different fields, for example see the books \cite{Corsini 2}, \cite{Davvaz book} and \cite{Vougiouklis}. Especially, the concept of hypervector space was introduced by Scafati-Tallini \cite{Tallini 1} in 1990 and has been investigated by herself, Ameri \cite{Ameri Dehghan 1}, Sedghi \cite{Sedghi Dehghan} and the author \cite{Dehghan lin func hvs,Dehghan Neutro HVS,Dehghan Ameri Ebrahimi}.

The concepts of fuzzy sets, soft sets, fuzzy soft sets and bipolar fuzzy soft sets have had mutual effects in algebraic hyperstructures (see the book \cite{Davvaz book 2}). For example, Ameri \cite{Ameri FHVS over VF} introduced a view of fuzzy hypervector spaces over valued fields in 2005 and studied some of their properties (\cite{Ameri Dehghan 2,Ameri Dehghan 3,Ameri Dehghan 4,Ameri Dehghan 5}). The author followed his work and studied some more properties of fuzzy hypervector spaces (\cite{Dehghan Sum SP on FHVS,Dehghan Various FQHVS,Dehghan bal abs,Dehghan afin conv,Dehghan Norouzi points}). Ranjbar \cite{Ranjbar} checked out some properties of soft hypervector spaces and fuzzy soft hypervector spaces. Xin \cite{Xin fuzzy soft BCK} applied the notion of intuitionistic fuzzy soft set to hyper BCK-algebras. The author \cite{Dehghan soft persian,Dehghan Nodehi soft} investigated some results in soft hypervector spaces. Norouzi \cite{Norouzi Ameri soft hmod} introduced some new directions on soft hypermodules and soft fuzzy hypermodules. Sarwar \cite{Sarwar} applied bipolar fuzzy soft sets to hypergraphs and elaborated various methods for the construction of bipolar fuzzy soft hypergraphs. Muhiuddin \cite{Muhiuddin} using the notion of bipolar-valued fuzzy soft set, investigated the concepts of bipolar-valued fuzzy soft hyper BCK-ideals.

Now in this paper, we apply the notion of bipolar fuzzy soft set in hypervector spaces and obtain some basic results which they will be a good basis for the future studies. More precisely, we define important operations on bipolar fuzzy soft sets, based on the operation on the hypervector space $V$. Then we define a bipolar fuzzy soft hypervector space, with some interesting examples, and investigate that the combinations of bipolar fuzzy soft hypervector spaces under the presented operation, are bipolar fuzzy soft hypervector spaces. Moreover, we check out the behaviour of bipolar fuzzy soft hypervector spaces under linear transformations.
\section{\protect\Large Preliminaries}

In this section we present some definitions and examples that we shall use in later.
\begin{definition}\cite{Zhang}
Let $X$ be a non-empty set. Then
\[A=\{(x,\mu^{+}_A(x),\mu^{-}_A(x)),\ x\in X\}\]
is called a bipolar fuzzy set in $X$, where $\mu^{+}_A:X\rightarrow [0,1]$ indicates that elements somewhat satisfy the property and $\mu^{-}_A:X\rightarrow [-1,0]$ indicates that elements somewhat satisfy the implicit counter-property.
\end{definition}
If $\mu^{+}_A(x)\neq 0$ and $\mu^{-}_A(x)=0$, it is the situation that $x$ has only positive satisfaction for $A$. If $\mu^{+}_A(x)=0$ and $\mu^{-}_A(x)\neq 0$, it means that $x$ does not satisfy the property of $A$, but it satisfies the counter-property of $A$.

For the sake of simplicity, we shall use the symbol $A=(\mu^{+}_A,\mu^{-}_A)$ or $A=(A^+,A^-)$, for the bipolar fuzzy set $A=\{(x,\mu^{+}_A(x),\mu^{-}_A(x)),\ x\in X\}$.

For example,
\[A=\{(mosquito,q,0),(dragon fly,0.4,0),(turtle,0,0),(snake,0,-1)\}\]
is a bipolar fuzzy set which represents the fuzzy concept frog's prey.
\begin{definition}\cite{Molodtsov}
Let $U$ be a universe set, $E$ be a set of parameters, $P(U)$ be the power set of $U$ and $A\subseteq E$. Then a pair $(F,A)$ is called a soft set over $U$, where $F$ is a mapping defined by $F:A\rightarrow P(U)$.
\end{definition}
\begin{definition}\cite{Abdullah}
Let $U$ be a universe set, $E$ be a set of parameters and $A\subseteq E$. Then a pair $(F,A)$ is said to be a bipolar fuzzy soft set over $U$, where $F$ is a mapping $F:A\rightarrow BF^U$ ($BF^U$ is the collection of all bipolar fuzzy sets over $U$), i.e.
\[\forall e\in A;\ F(e)=\{(x,\mu^{+}_{F(e)}(x),\mu^{-}_{F(e)}(x)),\ x\in U\}.\]
\end{definition}
For any $e\in A$, $F(e)$ is referred to as the set of $e$-approximate elements of the bipolar fuzzy soft set $(F,A)$, where $\mu^{+}_{F(e)}(x)$ denotes the degree of $x$ keeping the parameter $e$ and $\mu^{-}_{F(e)}(x)$ denotes the degree of $x$ keeping the non-parameter $e$. For simplicity, $\mu^{+}_{F(e)}(x)$ and $\mu^{-}_{F(e)}(x)$ are denoted by $F^{+}_{e}(x)$ and $F^{-}_{e}(x)$, respectively, and so
\[\forall e\in A;\ F(e)=F_{e}=\{(x,F^{+}_{e}(x),F^{-}_{e}(x)),\ x\in U\}.\]
\begin{definition}\label{D HVS}\cite{Tallini 1}
Let $K$ be a field, $(V,+)$ be an Abelian group and $P_{\ast}(V)$ be the set of all non-empty subsets of $V$. We define a hypervector space over $K$ to be the quadruplet $(V,+,\circ ,K)$, where $``\circ"$ is an external hyperoperation
\begin{equation*}
\circ :K\times V\longrightarrow P_{\ast }(V),
\end{equation*}
such that for all $a,b\in K$ and $x,y\in V$ the following conditions hold:
\begin{enumerate}
   \item[(H$_{1}$)] $a\circ (x+y)\subseteq a\circ x+a\circ y$, right distributive law,
   \item[(H$_{2}$)] $(a+b)\circ x\subseteq a\circ x+b\circ x$, left distributive law,
   \item[(H$_{3}$)] $a\circ (b\circ x)=(ab)\circ x$,
   \item[(H$_{4}$)] $a\circ (-x)=(-a)\circ x=-(a\circ x)$,
   \item[(H$_{5}$)] $x\in 1\circ x$,
\end{enumerate}
\noindent where in (H$_{1}$), $a\circ x+a\circ y=\{p+q:p\in a\circ x, q\in a\circ y\}$. Similarly it is in (H$_{2}$). Also in (H$_{3}$), $a\circ (b\circ x)=\underset{t\in b\circ x}{\bigcup}a\circ t$.

\noindent $V$ is called strongly right distributive, if we have equality in (H$_{1}$). In a similar way we define the strongly left distributive hypervector spaces.
\end{definition}
In the sequel of this paper, $V$ denotes a hypervector space over the field $K$, unless otherwise is specified.
\begin{example}\cite{Ameri Dehghan 3}\label{example hvs R3}
In classical vector space $(\mathbb{R}^{3},+,.,\mathbb{R})$ we define the external hyperoperation $\circ: \mathbb{R}\times \mathbb{R}^{3} \rightarrow P_{\ast}(\mathbb{R}^{3})$ by $a\circ(x_{0},y_{0},z_{0})=l$, where $``l"$ is a line with the parametric equations:
\begin{equation*}
l:\left\{
\begin{array}{l}
x=ax_{0}, \\
y=ay_{0}, \\
z=t.
\end{array}
\right.
\end{equation*}
Then $V=(\mathbb{R}^{3},+,\circ,\mathbb{R})$ is a strongly left distributive hypervector space over the field $\mathbb{R}$.
\end{example}
\begin{example}\label{example hvs Z4}
Let $K=\mathbb{Z}_2=\{0,1\}$ be the field of two numbers with the following operations:
\[
\begin{tabular}{ccc}
\begin{tabular}{c|c|c}
$+$ & $0$ & $1$ \\ \hline
$0$ & $0$ & $1$ \\ \hline
$1$ & $1$ & $0$ \\
\end{tabular}
& \ \ \ \ \ \ \ \ &
\begin{tabular}{c|c|c}
$\cdot $ & $0$ & $1$ \\ \hline
$0$ & $0$ & $0$ \\ \hline
$1$ & $0$ & $1$ \\
\end{tabular}
\end{tabular}
\]
Then $(\mathbb{Z}_4,+,\circ,\mathbb{Z}_2)$ is a hypervector space over the field $\mathbb{Z}_2$, where it is not strongly left or strongly right distributive, and the operation $``+:\mathbb{Z}_4\times \mathbb{Z}_4\rightarrow \mathbb{Z}_4"$ and the external hyperoperation $``\circ:\mathbb{Z}_2\times \mathbb{Z}_4\rightarrow P_{*}(\mathbb{Z}_4)"$ are defined as follow:
\[
\begin{array}{ccc}
\begin{tabular}{c|c|c|c|c}
$+$ & $0$ & $1$ & $2$ & $3$ \\ \hline
$0$ & $0$ & $1$ & $2$ & $3$ \\ \hline
$1$ & $1$ & $2$ & $3$ & $0$ \\ \hline
$2$ & $2$ & $3$ & $0$ & $1$ \\ \hline
$3$ & $3$ & $0$ & $1$ & $2$ \\
\end{tabular}
&\ \ \ \ \  &
\begin{tabular}{c|c|c|c|c}
$\circ $ & $0$ & $1$ & $2$ & $3$ \\ \hline
$0$ & $\{0,2\}$ & $\{0\}$ & $\{0\}$ & $\{0\}$ \\ \hline
$1$ & $\{0,2\}$ & $\{1,2,3\}$ & $\{0,2\}$ & $\{1,2,3\}$
\end{tabular}
\end{array}
\]
\end{example}
\section{\protect\Large Operations on Bipolar Fuzzy Soft Sets}
One of the first topics that is considered in the study of a set, from an algebraic point of view, is the definition of different operations on that set, in order to identify and study the created algebraic structure. In this section, this basic issue is discussed and some related properties are given.

Some operations have been defined on bipolar fuzzy soft sets by some authors. Here, at first we recall the definitions were defined by Akram \cite{Akram}, and then define new operations on bipolar fuzzy soft sets of hypervector spaces, based on their operation and external hyperoperation.

Let $(F,A)$ and $(G,B)$ be bipolar fuzzy soft sets over $U$. Then
\begin{enumerate}
  \item $(F,A)$ is called a bipolar fuzzy soft subset of $(G,B)$ and denoted by $(F,A)\sqsubseteq(G,B)$, if $A\subseteq B$ and for all $e\in A$, $F^{+}(e)\subseteq G^{+}(e)$ and $F^{-}(e)\supseteq G^{-}(e)$, i.e.
      \begin{quote}
        $\mu^{+}_{F(e)}(x)\leq\mu^{+}_{G(e)}(x)$ and $\mu^{-}_{F(e)}(x)\geq\mu^{-}_{G(e)}(x)$, for all $x\in U$.
      \end{quote}
  \item The intersection of $(F,A)$ and $(G,B)$ is denoted by $(F,A)\sqcap(G,B)$ and is defined as the bipolar fuzzy soft set $(F\sqcap G,A\cap B)$, where $(F\sqcap G)^{+}(e)=F^{+}(e)\cap G^{+}(e)$ and $(F\sqcap G)^{-}(e)=F^{-}(e)\cap G^{-}(e)$, for all $e\in A\cap B$, i.e.
      \begin{quote}
        $(F\sqcap G)^{+}_{e}(x)=F^{+}_{e}(x)\wedge G^{+}_{e}(x)$ and $(F\sqcap G)^{-}_{e}(x)=F^{-}_{e}(x)\vee G^{-}_{e}(x)$,
      \end{quote}
      for all $x\in U$.
  \item The extended intersection of $(F,A)$ and $(G,B)$ is denoted by $(F,A)\sqcap_\varepsilon(G,B)$ and is defined as the bipolar fuzzy soft set $(F\sqcap_\varepsilon G,A\cup B)$, where $(F\sqcap_\varepsilon G) (e)=F(e)$, for all $e\in A\setminus B$, $(F\sqcap_\varepsilon G) (e)=G(e)$, for all $e\in B\setminus A$, and $(F\sqcap_\varepsilon G)^+ (e)=F^+(e)\cap G^+(e)$, $(F\sqcap_\varepsilon G)^- (e)=F^-(e)\cap G^-(e)$, for all $e\in A\cap B$.
  \item The union of $(F,A)$ and $(G,B)$ is denoted by $(F,A)\sqcup(G,B)$ and is defined as the bipolar fuzzy soft set $(F\sqcup G,A\cup B)$, where $(F\sqcup G)(e)=F(e)$, for all $e\in A\setminus B$, $(F\sqcup G)(e)=G(e)$, for all $e\in B\setminus A$, and $(F\sqcup G)^+ (e)=F^+(e)\cup G^+(e)$, $(F\sqcup G)^- (e)= F^-(e)\cup G^-(e)$, for all $e\in A\cap B$, i.e.
      \begin{quote}
        $(F\sqcup G)^{+}_{e}(x)=F^{+}_{e}(x)\vee G^{+}_{e}(x)$ and $(F\sqcup G)^{-}_{e}(x)=F^{-}_{e}(x)\wedge G^{-}_{e}(x)$,
      \end{quote}
      for all $x\in U$.
  \item The restricted union of $(F,A)$ and $(G,B)$ is denoted by $(F,A)\sqcup_R(G,B)$ and is defined as the bipolar fuzzy soft set $(F\sqcup_R G,A\cap B)$, where $(F\sqcup_R G)^+ (e)=F^+(e)\cup G^+(e)$ and $(F\sqcup_R G)^- (e)=F^-(e)\cup G^-(e)$, for all $e\in A\cap B$.
  \item $(F,A)\ AND\ (G,B)$ is denoted by $(F,A)\wedge(G,B)$ and is defined as the bipolar fuzzy soft set $(F\wedge G,A\times B)$, where $(F\wedge G)^+(e_1,e_2)=F^+(e_1)\cap G^+(e_2)$ and $(F\wedge G)^- (e_1,e_2)=F^-(e_1)\cap G^-(e_2)$, for all $(e_1,e_2)\in A\times B$, i.e.
      \begin{equation*}
        (F\wedge G)^+_{(e_1,e_2)}(x)=F^+_{e_1}(x)\wedge G^+_{e_2}(x),
      \end{equation*}
      \begin{equation*}
        (F\wedge G)^-_{(e_1,e_2)}(x)= F^-_{e_1}(x)\vee G^-_{e_2}(x),
      \end{equation*}
      for all $x\in U$.
  \item $(F,A)\ OR\ (G,B)$ is denoted by $(F,A)\vee(G,B)$ and is defined as the bipolar fuzzy soft set $(F\vee G,A\times B)$, where $(F\vee G)^+(e_1,e_2)=F^+(e_1)\cup G^+(e_2)$ and $(F\wedge G)^- (e_1,e_2)=F^-(e_1) \cup G^-(e_2)$, for all $(e_1,e_2)\in A\times B$, i.e.
      \begin{equation*}
        (F\vee G)^+_{(e_1,e_2)}(x)=F^+_{e_1}(x)\vee G^+_{e_2}(x),
      \end{equation*}
      \begin{equation*}
        (F\vee G)^-_{(e_1,e_2)}(x)= F^-_{e_1}(x)\wedge G^-_{e_2}(x),
      \end{equation*}
      for all $x\in U$.
\end{enumerate}

Corresponding to every (hyper)operation over an algebraic (hyper)structure $U$, one can define an operation over bipolar fuzzy soft sets over $U$. In this paper, we study this idea in hypervector spaces.
\begin{definition}\label{def sum product bf soft set}
Let $(F,A)$ and $(G,B)$ be bipolar fuzzy soft sets over a hypervector space $V=(V,+,\circ,K)$ and $a\in K$. Then
\begin{enumerate}
  \item The sum of $(F,A)$ and $(G,B)$ is denoted by $(F,A)+(G,B)$ and is defined as the bipolar fuzzy soft set $(F+G,A\cap B)$, where
      \[(F+G)^+_e(x)=\bigvee_{x=y+z}(F^+_e(y)\wedge G^+_e(z)),\]
      \[(F+G)^-_e(x)=\bigwedge_{x=y+z}(F^-_e(y)\vee G^-_e(z)),\]
      for all $e\in A\cap B$ and $x\in V$.
  \item The extended sum of $(F,A)$ and $(G,B)$ is denoted by $(F,A)+_\varepsilon(G,B)$ and is defined as the bipolar fuzzy soft set $(F+_\varepsilon G,A\cup B)$, where
      \[
      (F+_{\varepsilon }G)_{e}^{+}(x)=\left\{
      \begin{array}{cl}
        F_{e}^{+}(x) & x\in A\setminus B, \\
        G_{e}^{+}(x) & x\in B\setminus A, \\
        (F+G)_{e}^{+}(x) & x\in A\cap B,
      \end{array}
      \right.
      \]
      and
      \[
      (F+_{\varepsilon }G)_{e}^{-}(x)=\left\{
      \begin{array}{cl}
        F_{e}^{-}(x) & x\in A\setminus B, \\
        G_{e}^{-}(x) & x\in B\setminus A, \\
        (F+G)_{e}^{-}(x) & x\in A\cap B,
      \end{array}
      \right.
      \]
      for all $e\in A\cup B$ and $x\in V$.
  \item The scalar product $a\circ(F,A)$ is defined as the bipolar fuzzy soft set $(a\circ F,A)$, where
    \[
    (a\circ F)_{e}^{+}(x)=\left\{
    \begin{array}{cl}
      \bigvee\limits_{x\in a\circ t}F_{e}^{+}(t) & \exists t\in V,x\in a\circ t, \\
      0 & otherwise,
    \end{array}
    \right.
    \]
  and
    \[
    (a\circ F)_{e}^{-}(x)=\left\{
    \begin{array}{cl}
    \bigwedge\limits_{x\in a\circ t}F_{e}^{-}(t) & \exists t\in V,x\in a\circ t, \\
    0 & otherwise.
    \end{array}
    \right.
    \]
\end{enumerate}
\end{definition}
\begin{proposition}
Let $(F,A)$ and $(G,B)$ be bipolar fuzzy soft sets of hypervector space $V=(V,+,\circ,K)$. Then for all $x,y\in V$ and $e\in A\cap B$, $(F+G)_e^+(x+y)\geq F_e^+(x)\wedge G_e^+(y)$ and $(F+G)_e^-(x+y)\leq F_e^-(x)\vee G_e^-(y)$.
\end{proposition}
\begin{proof}
By Definition \ref{def sum product bf soft set}, it follows that:
\[(F+G)^+_e(x+y)=\bigvee_{x+y=t_1+t_2}(F^+_e(t_1)\wedge G^+_e(t_2))\geq F_e^+(x)\wedge G_e^+(y),\]
and
\[(F+G)^-_e(x+y)=\bigwedge_{x+y=s_1+s_2}(F^-_e(s_1)\vee G^-_e(s_2))\leq F_e^-(x)\vee G_e^-(y).\]
\end{proof}
\begin{proposition}
Let $(F,A)$ be a bipolar fuzzy soft set of hypervector space $V=(V,+,\circ,K)$. Then $(F,A)\sqsubseteq 1\circ (F,A)$ and $-(F,A)\sqsubseteq (-1)\circ (F,A)$, where $(-F)_e^+(x)=F_e^+(-x)$ and $(-F)_e^-(x)=F_e^-(-x)$, for all $e\in A$, $x\in V$.
\end{proposition}
\begin{proof}
Let $e\in A$ and $x\in V$. Then by Definition \ref{D HVS}, $x\in 1\circ x$ and $-x\in 1\circ (-x)=(-1)\circ x$. Thus by Definition \ref{def sum product bf soft set}, it follows that:
\begin{equation*}
(1\circ F)_{e}^{+}(x)=\bigvee\limits_{x\in 1\circ t}F_{e}^{+}(t)\geq F_{e}^{+}(x),\ (1\circ F)_{e}^{-}(x)=\bigwedge\limits_{x\in 1\circ t}F_{e}^{-}(t)\leq F_{e}^{-}(x),
\end{equation*}
\begin{equation*}
((-1)\circ F)_{e}^{+}(x)=\bigvee\limits_{x\in (-1)\circ t}F_{e}^{+}(t)\geq F_{e}^{+}(-x)=(-F)_{e}^{+}(x),
\end{equation*}
and
\begin{equation*}
((-1)\circ F)_{e}^{-}(x)=\bigwedge\limits_{x\in (-1)\circ t}F_{e}^{-}(t)\leq F_{e}^{-}(-x)=(-F)_{e}^{-}(x).
\end{equation*}
\end{proof}
\begin{proposition}
Let $\{(F_i,A)\}_{i\in I}$, $\{(G_j,A)\}_{j\in J}$ be families of bipolar fuzzy soft sets over the hypervector space $V$ and $a\in K$. Then the followings hold:
\begin{enumerate}
  \item $\left(\bigsqcup\limits_{i\in I}(F_i,A)\right)+\left(\bigsqcup\limits_{j\in J}(G_j,A)\right)=\bigsqcup \limits_{i\in I,j\in J}((F_i,A)+ (G_j,A))$.
  \item $a\circ \left(\bigsqcup\limits_{i\in I}(F_i,A)\right)=\bigsqcup\limits_{i\in I}(a\circ (F_i,A))$.
\end{enumerate}
\end{proposition}
\begin{proof} 1) Let $e\in A$ and $x\in V$. Then
\begin{eqnarray*}
&&\left( \left( \bigsqcup\limits_{i\in I}(F_{i},A)\right) +\left(
\bigsqcup\limits_{j\in J}(G_{j},A)\right) \right) _{e}^{+}(x) \\
&=&\bigvee\limits_{x=y+z}\left( \left( \bigsqcup\limits_{i\in
I}(F_{i},A)\right) _{e}^{+}(y)\wedge \left( \bigsqcup\limits_{j\in
J}(G_{j},A)\right) _{e}^{+}(z)\right)  \\
&=&\bigvee\limits_{x=y+z}\left( \left( \bigvee\limits_{i\in
I}(F_{i})_{e}^{+}(y)\right) \wedge \left( \bigvee\limits_{j\in
J}(G_{j})_{e}^{+}(z)\right) \right)  \\
&=&\bigvee\limits_{x=y+z}\bigvee\limits_{i\in I,j\in J}\left(
(F_{i})_{e}^{+}(y)\wedge (G_{j})_{e}^{+}(z)\right)  \\
&=&\bigvee\limits_{i\in I,j\in J}\bigvee\limits_{x=y+z}\left(
(F_{i})_{e}^{+}(y)\wedge (G_{j})_{e}^{+}(z)\right)  \\
&=&\bigvee\limits_{i\in I,j\in J}((F_{i},A)+(G_{j},A))_{e}^{+}(x) \\
&=&\left( \bigsqcup\limits_{i\in I,j\in J}((F_{i},A)+(G_{j},A))\right)
_{e}^{+}(x),
\end{eqnarray*}
  and
 \begin{eqnarray*}
&&\left( \left( \bigsqcup\limits_{i\in I}(F_{i},A)\right) +\left(
\bigsqcup\limits_{j\in J}(G_{j},A)\right) \right) _{e}^{-}(x) \\
&=&\bigwedge\limits_{x=y+z}\left( \left( \bigsqcup\limits_{i\in
I}(F_{i},A)\right) _{e}^{-}(y)\vee \left( \bigsqcup\limits_{j\in
J}(G_{j},A)\right) _{e}^{-}(z)\right)  \\
&=&\bigwedge\limits_{x=y+z}\left( \left( \bigvee\limits_{i\in
I}(F_{i})_{e}^{-}(y)\right) \vee \left( \bigvee\limits_{j\in
J}(G_{j})_{e}^{-}(z)\right) \right)  \\
&=&\bigwedge\limits_{x=y+z}\bigwedge\limits_{i\in I,j\in J}\left(
(F_{i})_{e}^{-}(y)\vee (G_{j})_{e}^{-}(z)\right)  \\
&=&\bigwedge\limits_{i\in I,j\in J}\bigwedge\limits_{x=y+z}\left(
(F_{i})_{e}^{-}(y)\vee (G_{j})_{e}^{-}(z)\right)  \\
&=&\bigwedge\limits_{i\in I,j\in J}((F_{i},A)+(G_{j},A))_{e}^{-}(x) \\
&=&\left( \bigsqcup\limits_{i\in I,j\in J}((F_{i},A)+(G_{j},A))\right)
_{e}^{-}(x).
\end{eqnarray*}
2) Let $e\in A$ and $x\in V$. Then
  \begin{eqnarray*}
  \left( a\circ \left( \bigsqcup\limits_{i\in I}(F_{i},A)\right) \right)_{e}^{+}(x) &=&\left\{
  \begin{array}{cl}
  \bigvee\limits_{x\in a\circ t}\left( \bigsqcup\limits_{i\in I}(F_{i},A)\right) _{e}^{+}(t) & \exists t\in V,x\in a\circ t, \\
  0 & otherwise,
  \end{array}
  \right.  \\
  &=&\left\{
  \begin{array}{cl}
  \bigvee\limits_{x\in a\circ t}\left( \bigvee\limits_{i\in I}(F_{i})_{e}^{+}(t)\right)  & \exists t\in V,x\in a\circ t, \\
  0 & otherwise,
  \end{array}
  \right.  \\
  &=&\left\{
  \begin{array}{cl}
  \bigvee\limits_{i\in I}\left( \bigvee\limits_{x\in a\circ t}(F_{i})_{e}^{+}(t)\right)  & \exists t\in V,x\in a\circ t, \\
  0 & otherwise,
  \end{array}
  \right.  \\
  &=&\bigvee\limits_{i\in I}(a\circ (F_{i},A))_{e}^{+}(x) \\
  &=&\left( \bigsqcup\limits_{i\in I}(a\circ (F_{i},A))\right) _{e}^{+}(x),
  \end{eqnarray*}
  and
  \begin{eqnarray*}
  \left( a\circ \left( \bigsqcup\limits_{i\in I}(F_{i},A)\right) \right)_{e}^{-}(x) &=&\left\{
  \begin{array}{cl}
  \bigwedge\limits_{x\in a\circ t}\left( \bigsqcup\limits_{i\in I}(F_{i},A)\right) _{e}^{-}(t) & \exists t\in V,x\in a\circ t, \\
  0 & otherwise,
  \end{array}
  \right.  \\
  &=&\left\{
  \begin{array}{cl}
  \bigwedge\limits_{x\in a\circ t}\left( \bigwedge\limits_{i\in I}(F_{i})_{e}^{-}(t)\right)  & \exists t\in V,x\in a\circ t, \\
  0 & otherwise,
  \end{array}
  \right.  \\
  &=&\left\{
  \begin{array}{cl}
  \bigwedge\limits_{i\in I}\left( \bigwedge\limits_{x\in a\circ t}(F_{i})_{e}^{-}(t)\right)  & \exists t\in V,x\in a\circ t, \\
  0 & otherwise,
  \end{array}
  \right.  \\
  &=&\bigwedge\limits_{i\in I}(a\circ (F_{i},A))_{e}^{-}(x) \\
  &=&\left( \bigsqcup\limits_{i\in I}(a\circ (F_{i},A))\right) _{e}^{-}(x).
  \end{eqnarray*}
\end{proof}
\section{\protect\Large Bipolar Fuzzy Soft Hypervector Spaces}
In this section, the notion of bipolar fuzzy soft hypervector space is defined, based on especial bipolar fuzzy subsets of a hypervector space $V=(V,+,\circ,K)$, supported by some non-trivial examples. Moreover, it will be shown that intersection, extended intersection, union, restricted union, AND, sum, extended sum and scalar product of bipolar fuzzy soft hypervector spaces are bipolar fuzzy soft hypervector spaces, too. In fact, by the mentioned operations, new bipolar fuzzy soft hypervector spaces are constructed.
\begin{definition}\label{def bfshs}
Let $V$ be a hypervector space over the field $K$. Then a bipolar fuzzy set $A=(A^+,A^-)$ in $V$ is called a bipolar fuzzy subhyperspace of $V$, if for all $x,y\in V$ and $a\in K$ the followings hold:
\begin{enumerate}
  \item $A^{+}(x-y)\geq A^{+}(x)\wedge A^{+}(y)$, $A^{-}(x-y)\leq A^{-}(x)\vee A^{-}(y)$,
  \item $\bigwedge\limits_{t\in a\circ x}A^{+}(t)\geq A^{+}(x)$, $\bigvee\limits_{t\in a\circ x}A^{-}(t)\leq A^{-}(x)$.
\end{enumerate}
\end{definition}
\begin{example}\label{example bfshs R3}
Consider the hypervector space $V=(\mathbb{R}^3,+,\circ,\mathbb{R})$ in Example \ref{example hvs R3}. Define a bipolar fuzzy set $A=(A^+,A^-)$ in $V$, where $``A^+:\mathbb{R}^3\rightarrow [0,1]"$ and $``A^-:\mathbb{R}^3 \rightarrow [-1,0]"$ are given by
\[
A^{+}(x,y,z)=\left\{
\begin{array}{cl}
t_{1} & (x,y,z)\in \{0\}\times \{0\}\times\mathbb{R}, \\
t_{2} & (x,y,z)\in (\mathbb{R}\times \{0\}\times\mathbb{R})\setminus (\{0\}\times \{0\}\times\mathbb{R}), \\
t_{3} & otherwise,
\end{array}
\right.
\]
and
\[
A^{-}(x,y,z)=\left\{
\begin{array}{cl}
s_{1} & (x,y,z)\in \{0\}\times \{0\}\times\mathbb{R}, \\
s_{2} & (x,y,z)\in (\mathbb{R}\times \{0\}\times\mathbb{R})\setminus (\{0\}\times \{0\}\times\mathbb{R}), \\
s_{3} & otherwise,
\end{array}
\right.
\]
for some $-1\leq s_1<s_2<s_3\leq 0\leq t_3<t_2< t_1\leq 1$. Then $A=(A^+,A^-)$ is a bipolar fuzzy subhyperspace of $\mathbb{R}^3$.
\end{example}
\begin{example}\label{example bfshs Z4}
Consider the hypervector space $V=(\mathbb{Z}_4,+,\circ,\mathbb{Z}_2)$ in Example \ref{example hvs Z4}. Define a bipolar fuzzy set $B=(B^+,B^-)$ in $V$, where $``B^+:\mathbb{Z}_4\rightarrow [0,1]"$ and $``B^-:\mathbb{Z}_4 \rightarrow [-1,0]"$ are given by
\[
\begin{array}{ccc}
B^{+}(x)=\left\{
\begin{array}{cc}
t_{1} & x\in\{0,2\} \\
t_{2} & x\in\{1,3\}
\end{array}
\right.
& \ \ \ \ \ &
B^{-}(x)=\left\{
\begin{array}{cc}
s_{1} & x\in\{0,2\} \\
s_{2} & x\in\{1,3\}
\end{array}
\right.
\end{array}
\]
for some $-1\leq s_1<s_2\leq 0\leq t_2< t_1\leq 1$. Then $B=(B^+,B^-)$ is a bipolar fuzzy subhyperspace of $\mathbb{Z}_4$.
\end{example}
\begin{definition}\label{def bf soft hvs}
Let $(F,A)$ be a bipolar fuzzy soft set of a hypervector space $V=(V,+,\circ,K)$. Then $(F,A)$ is said to be a bipolar fuzzy soft hypervector space of $V$, if $F(e)$ is a bipolar fuzzy subhyperspace of $V$, for all $e\in A$, i.e.
\begin{enumerate}
  \item $F^{+}_e(x-y)\geq F^{+}_e(x)\wedge F^{+}_e(y)$, $F^{-}_e(x-y)\leq F^{-}_e(x)\vee F^{-}_e(y)$,
  \item $\bigwedge\limits_{t\in a\circ x}F^{+}_e(t)\geq F^{+}_e(x)$, $\bigvee\limits_{t\in a\circ x}F^{-}_e(t)\leq F^{-}_e(x)$.
\end{enumerate}
\end{definition}
\begin{example}\label{example bf soft hvs R3}
Consider the hypervector space $V=(\mathbb{R}^{3},+,\circ,\mathbb{R})$ in Example \ref{example hvs R3}. Suppose $A=\{a,b\}$ be a set of parameters. Then $(F,A)$ is a bipolar fuzzy soft hypervector space of $V$, where $``F^+_a,F^+_b:\mathbb{R}^{3}\rightarrow [0,1]"$ and $``F^-_a,F^-_b:\mathbb{R}^{3}\rightarrow [-1,0]"$ are given by the followings:
\[
F_{a}^{+}(x,y,z)=\left\{
\begin{array}{cl}
0.7 & (x,y,z)\in \{0\}\times \{0\}\times \mathbb{R}, \\
0.3 & (x,y,z)\in (\mathbb{R}\times \{0\}\times \mathbb{R})\setminus (\{0\}\times \{0\}\times \mathbb{R}), \\
0 & otherwise,
\end{array}
\right.
\]

\[
F_{a}^{-}(x,y,z)=\left\{
\begin{array}{cl}
-0.8 & (x,y,z)\in \{0\}\times \{0\}\times \mathbb{R}, \\
-0.4 & (x,y,z)\in (\mathbb{R}\times \{0\}\times \mathbb{R})\setminus (\{0\}\times \{0\}\times \mathbb{R}), \\
-0.2 & otherwise,
\end{array}
\right.
\]

\[
F_{b}^{+}(x,y,z)=\left\{
\begin{array}{cl}
0.9 & (x,y,z)\in \{0\}\times \{0\}\times \mathbb{R}, \\
0.4 & (x,y,z)\in (\mathbb{R}\times \{0\}\times \mathbb{R})\setminus (\{0\}\times \{0\}\times \mathbb{R}), \\
0.1 & otherwise,
\end{array}
\right.
\]

\[
F_{b}^{-}(x,y,z)=\left\{
\begin{array}{cl}
-0.6 & (x,y,z)\in \{0\}\times \{0\}\times \mathbb{R}, \\
-0.5 & (x,y,z)\in (\mathbb{R}\times \{0\}\times \mathbb{R})\setminus (\{0\}\times \{0\}\times \mathbb{R}), \\
-0.1 & otherwise.
\end{array}
\right.
\]
\end{example}
\begin{example}\label{example bf soft hvs Z4}
Consider the hypervector space $V=(\mathbb{Z}_4,+,\circ,\mathbb{Z}_2)$ in Example \ref{example hvs Z4}. Suppose $A=\{c,d,e\}$ be a set of parameters. Then $(F,A)$ is a bipolar fuzzy soft hypervector space of $V$, where $``F^+_c,F^+_d,F^+_e:\mathbb{Z}_4\rightarrow [0,1]"$ and $``F^-_c,F^-_d,F^-_e:\mathbb{Z}_4\rightarrow [-1,0]"$ are given by the followings:
\[
\begin{array}{ccc}
F_{c}^{+}(x)=\left\{
\begin{array}{cc}
0.5 & x\in \{0,2\} \\
0.3 & x\in \{1,3\}
\end{array}
\right.
& \ \ \ \ \  &
F_{c}^{-}(x)=\left\{
\begin{array}{cc}
-0.4 & x\in \{0,2\} \\
-0.2 & x\in \{1,3\}
\end{array}
\right.
\end{array}
\]

\[
\begin{array}{ccc}
F_{d}^{+}(x)=\left\{
\begin{array}{cc}
0.7 & x\in \{0,2\} \\
0.2 & x\in \{1,3\}
\end{array}
\right.
& \ \ \ \ \  &
F_{d}^{-}(x)=\left\{
\begin{array}{cc}
-0.6 & x\in \{0,2\} \\
-0.3 & x\in \{1,3\}
\end{array}
\right.
\end{array}
\]

\[
\begin{array}{ccc}
F_{e}^{+}(x)=\left\{
\begin{array}{cc}
0.8 & x\in \{0,2\} \\
0.4 & x\in \{1,3\}
\end{array}
\right.
& \ \ \ \ \  &
F_{e}^{-}(x)=\left\{
\begin{array}{cc}
-0.7 & x\in \{0,2\} \\
-0.5 & x\in \{1,3\}
\end{array}
\right.
\end{array}
\]
\end{example}
\begin{proposition}
Let $(F,A)$ and $(G,B)$ be bipolar fuzzy soft hypervector spaces of $V=(V,+,\circ,K)$. Then
\begin{enumerate}
  \item $(F,A)\sqcap(G,B)$ is a bipolar fuzzy soft hypervector space of $V$.
  \item $(F,A)\sqcap_\varepsilon(G,B)$ is a bipolar fuzzy soft hypervector space of $V$.
  \item If $A\cap B=\emptyset$, then $(F,A)\sqcup(G,B)$ is a bipolar fuzzy soft hypervector space of $V$.
  \item $(F,A)\sqcup_R(G,B)$ is a bipolar fuzzy soft hypervector space of $V$.
  \item $(F,A)\wedge(G,B)$ is a bipolar fuzzy soft hypervector space of $V$.
\end{enumerate}
\end{proposition}
\begin{proof}
We check the conditions of Definition \ref{def bf soft hvs}, for every item.
\begin{enumerate}
\item Let $e\in A\cap B$, $x,y\in V$ and $a\in K$. Then
\begin{eqnarray*}
(F\sqcap G)_{e}^{+}(x-y) &=&F_{e}^{+}(x-y)\wedge G_{e}^{+}(x-y) \\
&\geq &\left( F_{e}^{+}(x)\wedge F_{e}^{+}(y)\right) \wedge \left(G_{e}^{+}(x)\wedge G_{e}^{+}(y)\right)  \\
&=&\left( F_{e}^{+}(x)\wedge G_{e}^{+}(x)\right) \wedge \left(F_{e}^{+}(y)\wedge G_{e}^{+}(y)\right)  \\
&=&(F\sqcap G)_{e}^{+}(x)\wedge (F\sqcap G)_{e}^{+}(y),
\end{eqnarray*}
\begin{eqnarray*}
(F\sqcap G)_{e}^{-}(x-y) &=&F_{e}^{-}(x-y)\vee G_{e}^{-}(x-y) \\
&\leq &\left( F_{e}^{-}(x)\vee F_{e}^{-}(y)\right) \vee \left(G_{e}^{-}(x)\vee G_{e}^{-}(y)\right)  \\
&=&\left( F_{e}^{-}(x)\vee G_{e}^{-}(x)\right) \vee \left( F_{e}^{-}(y)\vee G_{e}^{-}(y)\right)  \\
&=&(F\sqcap G)_{e}^{+}(x)\vee (F\sqcap G)_{e}^{+}(y),
\end{eqnarray*}
\begin{eqnarray*}
\bigwedge\limits_{t\in a\circ x}(F\sqcap G)_{e}^{+}(t)
&=&\bigwedge\limits_{t\in a\circ x}\left( F_{e}^{+}(t)\wedge G_{e}^{+}(t)\right)  \\
&=&\left( \bigwedge\limits_{t\in a\circ x}F_{e}^{+}(t)\right) \wedge \left(\bigwedge\limits_{t\in a\circ x}G_{e}^{+}(t)\right)  \\
&\geq &F_{e}^{+}(x)\wedge G_{e}^{+}(x) \\
&=&(F\sqcap G)_{e}^{+}(x),
\end{eqnarray*}
\begin{eqnarray*}
\bigvee\limits_{t\in a\circ x}(F\sqcap G)_{e}^{-}(t) &=&\bigvee\limits_{t\in a\circ x}\left( F_{e}^{-}(t)\vee G_{e}^{-}(t)\right)  \\
&=&\left( \bigvee\limits_{t\in a\circ x}F_{e}^{+}(t)\right) \vee \left(\bigvee\limits_{t\in a\circ x}G_{e}^{+}(t)\right)  \\
&\leq &F_{e}^{-}(x)\vee G_{e}^{-}(x) \\
&=&(F\sqcap G)_{e}^{-}(x).
\end{eqnarray*}

\item It is similar to the proof of part (1).
\item If $A\cap B=\emptyset$, then
\[
(F\sqcup G)_{e}=\left\{
\begin{array}{cc}
F_{e} & e\in A\setminus B \\
G_{e} & e\in B\setminus A
\end{array}
\right.
\]
and thus $(F,A)\sqcup(G,B)$ is a bipolar fuzzy soft hypervector space of $V$.
\item The proof is obvious.
\item Let $(e_1,e_2)\in A\times B$, $x,y\in V$ and $a\in K$. Then
\begin{eqnarray*}
(F\wedge G)_{(e_{1},e_{2})}^{+}(x-y) &=&F_{e_{1}}^{+}(x-y)\wedge G_{e_{2}}^{+}(x-y) \\
&\geq &\left( F_{e_{1}}^{+}(x)\wedge F_{e_{1}}^{+}(y)\right) \wedge \left(G_{e_{2}}^{+}(x)\wedge G_{e_{2}}^{+}(y)\right)  \\
&=&\left( F_{e_{1}}^{+}(x)\wedge G_{e_{2}}^{+}(x)\right) \wedge \left(F_{e_{1}}^{+}(y)\wedge G_{e_{2}}^{+}(y)\right)  \\
&=&(F\wedge G)_{(e_{1},e_{2})}^{+}(x)\wedge (F\wedge G)_{(e_{1},e_{2})}^{+}(y),
\end{eqnarray*}
\begin{eqnarray*}
(F\wedge G)_{(e_{1},e_{2})}^{-}(x-y) &=&F_{e_{1}}^{-}(x-y)\vee G_{e_{2}}^{-}(x-y) \\
&\leq &\left( F_{e_{1}}^{-}(x)\vee F_{e_{1}}^{-}(y)\right) \vee \left(G_{e_{2}}^{-}(x)\vee G_{e_{2}}^{-}(y)\right)  \\
&=&\left( F_{e_{1}}^{-}(x)\vee G_{e_{2}}^{-}(x)\right)\vee\left(F_{e_{1}}^{-}(y)\vee G_{e_{2}}^{-}(y)\right) \\
&=&(F\wedge G)_{(e_{1},e_{2})}^{-}(x)\vee (F\wedge G)_{(e_{1},e_{2})}^{-}(y),
\end{eqnarray*}
\begin{eqnarray*}
\bigwedge\limits_{t\in a\circ x}(F\wedge G)_{(e_{1},e_{2})}^{+}(t)
&=&\bigwedge\limits_{t\in a\circ x}\left( F_{e_{1}}^{+}(t)\wedge G_{e_{2}}^{+}(t)\right)  \\
&=&\left( \bigwedge\limits_{t\in a\circ x}F_{e_{1}}^{+}(t)\right) \wedge \left( \bigwedge\limits_{t\in a\circ x}G_{e_{2}}^{+}(t)\right)  \\
&\geq &F_{e_{1}}^{+}(x)\wedge G_{e_{2}}^{+}(x) \\
&=&(F\wedge G)_{(e_{1},e_{2})}^{+}(x),
\end{eqnarray*}
\begin{eqnarray*}
\bigvee\limits_{t\in a\circ x}(F\wedge G)_{(e_{1},e_{2})}^{-}(t)
&=&\bigvee\limits_{t\in a\circ x}\left( F_{e_{1}}^{-}(t)\vee G_{e_{2}}^{-}(t)\right)  \\
&=&\left( \bigvee\limits_{t\in a\circ x}F_{e_{1}}^{-}(t)\right) \vee \left(\bigvee\limits_{t\in a\circ x}G_{e_{2}}^{-}(t)\right)  \\
&\leq &F_{e_{1}}^{-}(x)\vee G_{e_{2}}^{-}(x) \\
&=&(F\wedge G)_{(e_{1},e_{2})}^{-}(x).
\end{eqnarray*}
\noindent Hence $(F,A)\wedge(G,B)$ is a bipolar fuzzy soft hypervector space of $V$.
\end{enumerate}
\end{proof}
\begin{proposition}\label{prop sum bf soft hvs}
Let $(F,A)$, $(G,B)$ be bipolar fuzzy soft hypervector spaces of $V=(V,+,\circ,K)$. Then $(F,A)+(G,B)$ is a bipolar fuzzy soft hypervector space of $V$.
\end{proposition}
\begin{proof}
Let $e\in A\cap B$, $x,y\in V$, $a\in K$,
\[A_x=\{(x_1,x_2)\in V\times V;\ x=x_1+x_2\},\]
\[A_y=\{(y_1,y_2)\in V\times V;\ y=y_1+y_2\}.\]
Hence

1) If $A_x=\emptyset$ or $A_y=\emptyset$, then obviously $(F+G)_{e}^{+}(x-y)\geq (F+G)_{e}^{+}(x)\wedge (F+G)_{e}^{+}(y)$. Otherwise, if $A_x\neq\emptyset$ or $A_y\neq\emptyset$, then there exist $x_1,x_2,y_1,y_2\in V$, such that $x=x_1+x_2$ and $y=y_1+y_2$. Thus
\[A_{x-y}=\{(t_1,t_2)\in V\times V;\ x-y=t_1+t_2\}\neq\emptyset,\]
and so
\begin{eqnarray*}
(F+G)_{e}^{+}(x-y) &=&\bigvee\limits_{x-y=t_{1}+t_{2}}\left(F_{e}^{+}(t_{1})\wedge G_{e}^{+}(t_{2})\right)  \\
&\geq &F_{e}^{+}(x_{1}-y_{1})\wedge G_{e}^{+}(x_{2}-y_{2}) \\
&\geq &\left( F_{e}^{+}(x_{1})\wedge F_{e}^{+}(y_{1})\right) \wedge \left(G_{e}^{+}(x_{2})\wedge G_{e}^{+}(y_{2})\right)  \\
&=&\left( F_{e}^{+}(x_{1})\wedge G_{e}^{+}(x_{2})\right) \wedge \left(F_{e}^{+}(y_{1})\wedge G_{e}^{+}(y_{2})\right).
\end{eqnarray*}
Hence
\begin{eqnarray*}
(F+G)_{e}^{+}(x-y) &\geq &\left( \bigvee\limits_{x=x_{1}+x_{2}}\left(F_{e}^{+}(x_{1})\wedge G_{e}^{+}(x_{2}) \right) \right)  \\
&&\wedge \left( \bigvee\limits_{y=y_{1}+y_{2}}\left( F_{e}^{+}(y_{1})\wedge G_{e}^{+}(y_{2})\right) \right)  \\
&=&(F+G)_{e}^{+}(x)\wedge (F+G)_{e}^{+}(y).
\end{eqnarray*}
For the negative part, if $A_x=\emptyset$ or $A_y=\emptyset$, then clearly $(F+G)_{e}^{-}(x-y)\leq (F+G)_{e}^{-}(x)\vee (F+G)_{e}^{-}(y)$. Otherwise, if $A_x\neq\emptyset$ or $A_y\neq\emptyset$, then similar to the previous part, it follows that
\begin{eqnarray*}
(F+G)_{e}^{-}(x-y) &=&\bigwedge\limits_{x-y=t_{1}+t_{2}}\left(F_{e}^{-}(t_{1})\vee G_{e}^{-}(t_{2})\right) \\
&\leq &F_{e}^{-}(x_{1}-y_{1})\vee G_{e}^{-}(x_{2}-y_{2}) \\
&\leq &\left(F_{e}^{-}(x_{1})\vee F_e^{-}(y_{1})\right)\vee\left(G_{e}^{-}(x_{2})\vee G_{e}^{-}(y_{2})\right)\\
&=&\left( F_{e}^{-}(x_{1})\vee G_{e}^{-}(x_{2})\right) \vee \left(F_{e}^{-}(y_{1})\vee G_{e}^{+}(y_{2})\right),
\end{eqnarray*}
and so
\begin{eqnarray*}
(F+G)_{e}^{-}(x-y) &\leq &\left( \bigwedge\limits_{x=x_{1}+x_{2}}\left(F_{e}^{-}(x_{1})\vee G_{e}^{-}(x_{2}) \right) \right)  \\
&&\vee \left( \bigwedge\limits_{y=y_{1}+y_{2}}\left( F_{e}^{-}(y_{1})\vee G_{e}^{-}(y_{2})\right) \right)  \\
&=&(F+G)_{e}^{-}(x)\vee (F+G)_{e}^{-}(y).
\end{eqnarray*}
2) If $A_x=\emptyset$, then it is clear that $\bigwedge\limits_{t\in a\circ x}(F+G)_{e}^{+}(t)\geq (F+G)_{e}^{+}(x)$. If $A_x\neq\emptyset$, and $x_1,x_2\in V$ such that $x=x_1+x_2$, then for all $t\in a\circ x$, $t\in a\circ(x_1+x_2)\subseteq a\circ x_1+a\circ x_2$. Thus there exist $\acute{t}_{1}, \acute{t}_{2}\in V$ such that $\acute{t}_{1}\in a\circ x_{1}$, $\acute{t}_{2}\in a\circ x_{2}$ and $t=\acute{t}_{1}+\acute{t}_{2}$. Hence
\begin{eqnarray*}
\bigvee\limits_{t=t_{1}+t_{2}}\left( F_{e}^{+}(t_{1})\wedge G_{e}^{+}(t_{2})\right)
&\geq & F_{e}^{+}(\acute{t}_{1})\wedge G_{e}^{+}(\acute{t}_{2}) \\
&\geq &\left( \bigwedge\limits_{r_{1}\in a\circ x_{1}}F_{e}^{+}(r_{1})\right) \wedge \left(\bigwedge\limits_ {r_{2}\in a\circ x_{2}}G_{e}^{+}(r_{2})\right) \\
&\geq &F_{e}^{+}(x_{1})\wedge G_{e}^{+}(x_{2}).
\end{eqnarray*}
It follows that,
\[
(F+G)_{e}^{+}(t)\geq \bigvee\limits_{x=x_{1}+x_{2}}\left(F_{e}^{+}(x_{1})\wedge G_{e}^{+}(x_{2})\right) =(F+G)_{e}^{+}(x),
\]
and so $\bigwedge\limits_{t\in a\circ x}(F+G)_{e}^{+}(t)\geq (F+G)_{e}^{+}(x)$.

Now, for the negative part, if $A_x=\emptyset$, then obviously $\bigvee\limits_{t\in a\circ x}(F+G)_{e}^{-}(t) \leq (F+G)_{e}^{-}(x)$, and if $A_x\neq\emptyset$ and $x_1,x_2\in V$ such that $x=x_1+x_2$, then for all $t\in a\circ x$, $t=\acute{t}_{1}+\acute{t}_{2}$, for some $\acute{t}_{1}\in a\circ x_{1}$ and $\acute{t}_{2}\in a\circ x_{2}$. Then
\begin{eqnarray*}
(F+G)_{e}^{-}(t) &=&\bigwedge\limits_{t=t_{1}+t_{2}}\left(F_{e}^{-}(t_{1})\vee G_{e}^{-}(t_{2})\right)  \\
&\leq &F_{e}^{-}(\acute{t}_{1})\vee G_{e}^{-}(\acute{t}_{2}) \\
&\leq &\left( \bigvee\limits_{s_{1}\in a\circ x_{1}}F_{e}^{-}(s_{1})\right)\vee \left( \bigvee\limits_{s_{2}\in a\circ x_{2}}G_{e}^{-}(s_{2})\right)  \\
&\leq &F_{e}^{-}(x_{1})\vee G_{e}^{-}(x_{2}).
\end{eqnarray*}
Thus $(F+G)_{e}^{-}(t)\leq \bigwedge\limits_{x=x_{1}+x_{2}}\left(F_{e}^{-}(x_{1})\vee G_{e}^{-}(x_{2})\right) =(F+G)_{e}^{-}(x)$, and so $\bigvee\limits_{t\in a\circ x}(F+G)_{e}^{-}(t)\leq (F+G)_{e}^{-}(x)$.

\noindent Therefor, the proof is completed.
\end{proof}
\begin{corollary}
Let $(F,A)$, $(G,B)$ be bipolar fuzzy soft hypervector spaces of $V=(V,+,\circ,K)$. Then $(F,A)+_\varepsilon (G,B)$ is a bipolar fuzzy soft hypervector space of $V$.
\end{corollary}
\begin{proof}
It is obvious, by Proposition \ref{prop sum bf soft hvs}.
\end{proof}
\begin{proposition}
Let $V=(V,+,\circ,K)$ be strongly right distributive and $(F,A)$ be a bipolar fuzzy soft hypervector space of $V$. Then $(a\circ F,A)$ is a bipolar fuzzy soft hypervector space of $V$, for all $a\in K$.
\end{proposition}
\begin{proof}
Let $e\in A$, $x,y\in V$ and $b\in K$. Then
\begin{enumerate}
  \item If there does not exist $t_1\in V$ such that $x\in a\circ t_1$, or there does not exist $t_2\in V$ such that $y\in a\circ t_2$, then it is clear that
      \[(a\circ F)_{e}^{+}(x-y)\geq (a\circ F)_{e}^{+}(x)\wedge (a\circ F)_{e}^{+}(y),\]
      and
      \[(a\circ F)_{e}^{-}(x-y)\leq (a\circ F)_{e}^{-}(x)\vee (a\circ F)_{e}^{-}(y).\]
      But, if there exist $t_{1},t_{2}\in V$, such that $x\in a\circ t_{1}$ and $y\in a\circ t_{2}$, then $x-y\in a\circ t_{1}-a\circ t_{2}=a\circ (t_{1}-t_{2})$, and so
      \begin{eqnarray*}
        (a\circ F)_{e}^{+}(x-y) &=&\bigvee\limits_{x-y\in a\circ t}F_{e}^{+}(t) \\
        &\geq &F_{e}^{+}(t_{1}-t_{2}) \\
        &\geq &F_{e}^{+}(t_{1})\wedge F_{e}^{+}(t_{2}),
      \end{eqnarray*}
      thus
      \begin{eqnarray*}
        (a\circ F)_{e}^{+}(x-y) &\geq &\left( \bigvee\limits_{x\in a\circ t_{1}}F_{e}^{+}(t_{1})\right) \wedge \left( \bigvee\limits_{y\in a\circ t_{2}}F_{e}^{+}(t_{2})\right)  \\
        &=&(a\circ F)_{e}^{+}(x)\wedge (a\circ F)_{e}^{+}(y).
      \end{eqnarray*}
      Also,
      \begin{eqnarray*}
        (a\circ F)_{e}^{-}(x-y) &=&\bigwedge\limits_{x-y\in a\circ t}F_{e}^{-}(t) \\
        &\leq &F_{e}^{-}(t_{1}-t_{2}) \\
        &\leq &F_{e}^{-}(t_{1})\vee F_{e}^{-}(t_{2}),
      \end{eqnarray*}
      hence
      \begin{eqnarray*}
        (a\circ F)_{e}^{-}(x-y) &\leq &\left( \bigwedge\limits_{x\in a\circ t_{1}}F_{e}^{-}(t_{1})\right) \vee \left( \bigwedge\limits_{y\in a\circ t_{2}}F_{e}^{-}(t_{2})\right)  \\
        &=&(a\circ F)_{e}^{-}(x)\vee (a\circ F)_{e}^{-}(y).
      \end{eqnarray*}
  \item If there does not exist $t\in V$ such that $x\in a\circ t$, then $(a\circ F)_{e}^{+}(x)=(a\circ F)_{e}^{-}(x)=0$ and clearly the result is obtained. Otherwise, if there exists $t\in V$, such that $x\in a\circ t$, then for all $s\in b\circ x$ and $r\in a\circ s$, it follows that
\[
F_{e}^{+}(r)\geq \bigwedge\limits_{l\in a\circ s}F_{e}^{+}(l)\geq F_{e}^{+}(s)\geq \bigwedge\limits_{k\in b\circ x}F_{e}^{+}(k)\geq F_{e}^{+}(x)\geq \bigwedge\limits_{p\in a\circ t}F_{e}^{+}(p)\geq F_{e}^{+}(t).
\]
Thus $F_{e}^{+}(r)\geq \bigvee\limits_{x\in a\circ t}F_{e}^{+}(t)$, and so $\bigvee\limits_{r\in a\circ s}F_{e}^{+}(r)\geq \bigvee\limits_{x\in a\circ t}F_{e}^{+}(t)$. Hence
\begin{eqnarray*}
\bigwedge\limits_{s\in b\circ x}(a\circ F)_{e}^{+}(s)
&=&\bigwedge\limits_{s\in b\circ x}\left( \bigvee\limits_{r\in a\circ s}F_{e}^{+}(r)\right)  \\
&\geq &\bigvee\limits_{x\in a\circ t}F_{e}^{+}(t) \\
&=&(a\circ F)_{e}^{+}(x).
\end{eqnarray*}
Moreover,
\[
F_{e}^{-}(r)\leq \bigvee\limits_{\acute{l}\in a\circ s}F_{e}^{-}(\acute{l})\leq F_{e}^{-}(s)\leq \bigvee\limits_{\acute{k}\in b\circ x}F_{e}^{-}(\acute{k})\leq F_{e}^{-}(x)\leq \bigvee\limits_{\acute{p}\in a\circ t}F_{e}^{-}(\acute{p})\leq F_{e}^{-}(t).
\]
Thus $F_{e}^{-}(r)\leq \bigwedge\limits_{x\in a\circ t}F_{e}^{-}(t)$, and so $\bigwedge\limits_{r\in a\circ s}F_{e}^{-}(r)\leq \bigwedge\limits_{x\in a\circ t}F_{e}^{-}(t)$. Hence
\begin{eqnarray*}
\bigvee\limits_{s\in b\circ x}(a\circ F)_{e}^{-}(s) &=&\bigvee\limits_{s\in b\circ x}\left( \bigwedge\limits_{r\in a\circ s}F_{e}^{-}(r)\right)  \\
&\leq &\bigwedge\limits_{x\in a\circ t}F_{e}^{-}(t) \\
&=&(a\circ F)_{e}^{-}(x).
\end{eqnarray*}
\end{enumerate}
Therefore, the proof is completed.
\end{proof}
\section{\protect\Large Bipolar Fuzzy Soft Hypervector Spaces under Linear Transformations}
In this section, by developing the extension principle of fuzzy sets to bipolar fuzzy soft sets, the behavior of bipolar fuzzy soft hypervector spaces under linear transformations is studied.
\begin{definition}
Let $(F,A)$ and $(G,B)$ be bipolar fuzzy soft sets over $X$ and $Y$, respectively. Then a pair $(\phi,f)$ is called a fuzzy soft function from $X$ to $Y$, where $\phi:X\rightarrow Y$ and $f:A\rightarrow B$ are functions. In this case, the image of $(F,A)$ under $(\phi,f)$ is the bipolar fuzzy soft set $(\phi,f)(F,A)=(\phi(F),f(A))$ of $Y$ defined by:
\begin{equation*}
\phi (F)_{u}^{+}(y)=\left\{
\begin{array}{cc}
\bigvee\limits_{x\in \phi ^{-1}(y)}\bigvee\limits_{e\in f^{-1}(u)}F_{e}^{+}(x) & \phi ^{-1}(y)\neq \phi, \\
0 & \phi ^{-1}(y)=\phi,
\end{array}
\right.
\end{equation*}
and
\begin{equation*}
\phi (F)_{u}^{-}(y)=\left\{
\begin{array}{cc}
\bigwedge\limits_{x\in \phi ^{-1}(y)}\bigwedge\limits_{e\in f^{-1}(u)}F_{e}^{-}(x) & \phi ^{-1}(y)\neq \phi, \\
0 & \phi ^{-1}(y)=\phi,
\end{array}
\right.
\end{equation*}
for all $u\in f(A)$, $y\in Y$.

\noindent Also, pre-image of $(G,B)$ under $(\phi,f)$ is the bipolar fuzzy soft set $(\phi,f)^{-1}(G,B)= (\phi^{-1}(G),A)$ of $X$, such that $\phi^{-1}(G)_e^+(x)=G_{f(e)}^+(\phi(x))$ and $\phi^{-1}(G)_e^-(x)= G_{f(e)} ^-(\phi(x))$, for all $e\in A$, $x\in X$.
\end{definition}
\begin{definition}
A fuzzy soft function $(T:V\rightarrow W,\ f:A\rightarrow B)$ between hypervector spaces $V$ and $W$ is said to be a bipolar fuzzy soft linear transformation if $T$ is a linear transformation, i.e. $T(x+y)=T(x)+T(y)$ and $T(a\circ x)\subseteq a\circ T(x)$, for all $x,y \in V$, $a\in K$. If $T(a\circ x)=a\circ T(x)$, then $(T,f)$ is said to be a bipolar fuzzy soft good transformation.
\end{definition}
\begin{theorem}
Let $(T:V\rightarrow W,\ f:A\rightarrow B)$ be a bipolar fuzzy soft good transformation. If $(F,A)$ is a bipolar fuzzy soft hypervector space of $V$, then $(T,f)(F,A)$ is a bipolar fuzzy soft hypervector space of $W$.
\end{theorem}
\begin{proof}
Let $u\in f(A)$, $\acute{x},\acute{y}\in W$, $a\in K$. If $T^{-1}(\acute{x})=\emptyset$ or $T^{-1}(\acute{y})= \emptyset$, then the proof is clear. If $T^{-1}(\acute{x})\neq\emptyset$ and $T^{-1}(\acute{y})\neq\emptyset$, then $T^{-1}(\acute{x}-\acute{y})\neq\emptyset$ and it follows that:
\begin{eqnarray*}
T(F)_{u}^{+}(\acute{x}-\acute{y})&=&\bigvee\limits_{T(t)=\acute{x}-\acute{y}}\bigvee\limits_{f(e)=u}F_{e}^{+}(t)\\
&\geq &\bigvee\limits_{T(x)=\acute{x},T(y)=\acute{y}}\bigvee\limits_{f(e)=u}F_{e}^{+}(x-y) \\
&\geq &\bigvee\limits_{T(x)=\acute{x},T(y)=\acute{y}}\bigvee\limits_{f(e)=u}\left( F_{e}^{+}(x)\wedge F_{e}^{+} (y)\right)  \\
&=&\left(\bigvee\limits_{T(x)=\acute{x}}\bigvee\limits_{f(e)=u}F_{e}^{+}(x)\right) \wedge \left(\bigvee\limits _{T(y)=\acute{y}}\bigvee\limits_{f(e)=u}F_{e}^{+}(y)\right)  \\
&=&T(F)_{u}^{+}(\acute{x})\wedge T(F)_{u}^{+}(\acute{y}),
\end{eqnarray*}
and%
\begin{eqnarray*}
T(F)_{u}^{-}(\acute{x}-\acute{y}) &=&\bigwedge\limits_{T(t)=\acute{x}-\acute{y}}\bigwedge\limits_{f(e)=u} F_{e}^{-}(t) \\
&\leq &\bigwedge\limits_{T(x)=\acute{x},T(y)=\acute{y}}\bigwedge\limits_{f(e)=u}F_{e}^{-}(x-y) \\
&\leq &\bigwedge\limits_{T(x)=\acute{x},T(y)=\acute{y}}\bigwedge\limits_{f(e)=u}\left( F_{e}^{-}(x)\vee F_{e}^{-}(y)\right)  \\
&=&\left( \bigwedge\limits_{T(x)=\acute{x}}\bigwedge\limits_{f(e)=u}F_{e}^{-}(x)\right) \vee \left(\bigwedge \limits_{T(y)=\acute{y}}\bigwedge\limits_{f(e)=u}F_{e}^{-}(y)\right)  \\
&=&T(F)_{u}^{-}(\acute{x})\vee T(F)_{u}^{-}(\acute{y}).
\end{eqnarray*}
Moreover, if $T^{-1}(\acute{x})=\emptyset$, then the proof is obvious. If $T^{-1}(\acute{x})\neq\emptyset$ and $x\in V$ such that $T(x)=\acute{x}$, then for any $\acute{t}\in a\circ \acute{x}$, $\acute{t}\in a\circ T(x)=T(a\circ x)$, so there exists $t\in a\circ x$, such that $\acute{t}=T(t)$. Thus $T^{-1}(\acute{t}) \neq\emptyset$ and it follows that:
\begin{eqnarray*}
T(F)_{u}^{+}(\acute{t}) &=&\bigvee\limits_{T(s)=\acute{t}}\bigvee\limits_{f(e)=u}F_{e}^{+}(s) \\
&=&\left( \bigvee\limits_{s\in a\circ x,T(s)=\acute{t}}\bigvee\limits_{f(e)=u}F_{e}^{+}(s)\right) \vee \left( \bigvee\limits_{s\in V\backslash a\circ x,T(s)=\acute{t}}\bigvee\limits_{f(e)=u}F_{e}^{+}(s)\right)  \\
&\geq &\bigvee\limits_{f(e)=u}\bigvee\limits_{s\in a\circ x,T(s)=\acute{t}}F_{e}^{+}(s) \\
&\geq &\bigvee\limits_{f(e)=u}\bigwedge\limits_{s\in a\circ x,T(s)=\acute{t}}F_{e}^{+}(s) \\
&\geq &\bigvee\limits_{f(e)=u}\bigwedge\limits_{s\in a\circ x}F_{e}^{+}(s) \\
&\geq &\bigvee\limits_{f(e)=u}F_{e}^{+}(x),
\end{eqnarray*}
Hence $T(F)_{u}^{+}(\acute{t})\geq \bigvee\limits_{T(x)=\acute{x}}\bigvee\limits_{f(e)=u}F_{e}^{+}(x) =T(F)_{u}^{+}(\acute{x})$, and so $\bigwedge\limits_{\acute{t}\in a\circ \acute{x}}T(F)_{u}^{+}(\acute{t})\geq
T(F)_{u}^{+}(\acute{x})$.

\noindent Moreover,
\begin{eqnarray*}
T(F)_{u}^{-}(\acute{t}) &=&\bigwedge\limits_{T(s)=\acute{t}}\bigwedge\limits_{f(e)=u}F_{e}^{-}(s) \\
&=&\left( \bigwedge\limits_{s\in a\circ x,T(s)=\acute{t}}\bigwedge\limits_{f(e)=u}F_{e}^{-}(s)\right) \wedge \left( \bigwedge\limits_{s\in V\backslash a\circ x,T(s)=\acute{t}}\bigwedge\limits_{f(e)=u}F_{e}^{-}(s)\right) \\
&\leq &\bigwedge\limits_{f(e)=u}\bigwedge\limits_{s\in a\circ x,T(s)=\acute{t}}F_{e}^{-}(s) \\
&\leq &\bigwedge\limits_{f(e)=u}\bigvee\limits_{s\in a\circ x,T(s)=\acute{t}}F_{e}^{-}(s) \\
&\leq &\bigwedge\limits_{f(e)=u}\bigvee\limits_{s\in a\circ x}F_{e}^{-}(s) \\
&\leq &\bigwedge\limits_{f(e)=u}F_{e}^{-}(x).
\end{eqnarray*}
Hence
\[T(F)_{u}^{-}(\acute{t})\leq \bigwedge\limits_{T(x)=\acute{x}}\bigwedge\limits_{f(e)=u}F_{e}^{-}(x)= T(F)_{u}^{-}(\acute{x}),\]
and so $\bigvee\limits_{\acute{t}\in a\circ \acute{x}}T(F)_{u}^{-}(\acute{t})\leq T(F)_{u}^{-}(\acute{x})$.

\noindent Therefore, by Definition \ref{def bf soft hvs}, $(T,f)(F,A)$ is a bipolar fuzzy soft hypervector space of $W$.
\end{proof}
\begin{theorem}
Let $(T:V\rightarrow W,\ f:A\rightarrow B)$ be a bipolar fuzzy soft linear transformation. If $(G,B)$ is a bipolar fuzzy soft hypervector space of $W$, then $(T,f)^{-1}(G,B)$ is a bipolar fuzzy soft hypervector space of $V$.
\end{theorem}
\begin{proof}
Let $e\in A$, $x,y\in V$, and $a\in K$. Then
\begin{eqnarray*}
T^{-1}(G)_{e}^{+}(x-y) &=&G_{f(e)}^{+}(T(x-y)) \\
&=&G_{f(e)}^{+}(T(x)-T(y)) \\
&\geq &G_{f(e)}^{+}(T(x))\wedge G_{f(e)}^{+}(T(y)) \\
&=&T^{-1}(G)_{e}^{+}(x)\wedge T^{-1}(G)_{e}^{+}(y),
\end{eqnarray*}
and
\begin{eqnarray*}
T^{-1}(G)_{e}^{-}(x-y) &=&G_{f(e)}^{-}(T(x-y)) \\
&=&G_{f(e)}^{-}(T(x)-T(y)) \\
&\leq &G_{f(e)}^{-}(T(x))\vee G_{f(e)}^{-}(T(y)) \\
&=&T^{-1}(G)_{e}^{-}(x)\vee T^{-1}(G)_{e}^{-}(y).
\end{eqnarray*}
Also,
\begin{eqnarray*}
\bigwedge\limits_{t\in a\circ x}T^{-1}(G)_{e}^{+}(t)&=&\bigwedge\limits_{t\in a\circ x}G_{f(e)}^{+}(T(t)) \\
&\geq &\bigwedge\limits_{\acute{t}\in a\circ T(x)}G_{f(e)}^{+}(\acute{t}) \\
&\geq &G_{f(e)}^{+}(T(x)) \\
&=&T^{-1}(G)_{e}^{+}(x),
\end{eqnarray*}
and
\begin{eqnarray*}
\bigvee\limits_{t\in a\circ x}T^{-1}(G)_{e}^{-}(t) &=&\bigvee\limits_{t\in a\circ x}G_{f(e)}^{-}(T(t)) \\
&\leq &\bigvee\limits_{\acute{t}\in a\circ T(x)}G_{f(e)}^{-}(\acute{t}) \\
&\leq &G_{f(e)}^{-}(T(x)) \\
&=&T^{-1}(G)_{e}^{-}(x).
\end{eqnarray*}
Therefore, by Definition \ref{def bf soft hvs}, $(T,f)^{-1}(G,B)$ is a bipolar fuzzy soft hypervector space of $V$.
\end{proof}
\section{\protect\Large Conclusion}
Various methods are known for mathematical modeling of imprecise phenomena. Two of these methods are the use of soft sets and bipolar fuzzy sets. Combining these two concepts with each other, i.e. bipolar fuzzy soft sets, leads to more accurate modeling of the discussed concepts. By defining different operations on the discussed set, an algebraic structure is obtained. In this article, these issues are mixed together and the application of bipolar fuzzy soft sets modeling method on the algebraic structure of hypervector space is studied. In this regard, while introducing the algebraic structure of bipolar fuzzy soft hypervector space, which is supported with interesting examples, some of it's basic features have been investigated and the basis for further studies has been prepared. With this information, the following topics can be studied in the future:

- Equivalent conditions of fuzzy bipolar soft hypervector space,

- Bipolar fuzzy soft hypervector spaces generated by a bipolar fuzzy soft set,

- Normal bipolar fuzzy soft hypervector spaces,

- Cosets of bipolar fuzzy soft hypervector spaces,

- Bipolar fuzzy soft sets on quotient hypervector spaces,

- Finding the applications of the introduced structure in decision making,

- Application of bipolar fuzzy soft sets over other algebraic structures/hyperstructures.

\vspace{0.7cm}
\noindent {\small\emph{Department of Mathematics, Faculty of Basic Sciences, University of Bojnord, Bojnord, Iran}}

\noindent \emph{dehghan@ub.ac.ir}

\end{document}